\documentclass[12pt]{article}

\usepackage{amssymb,latexsym}
\usepackage{amsfonts}
\usepackage{amsmath,amsthm,graphicx}

\newtheorem{thm}{Theorem}[section]
\newtheorem{lem}[thm]{Lemma}

\newtheorem{defn}[thm]{Definition}

\newtheorem{cor}[thm]{Corollary}

\newcommand{\integers}{{\mathbb Z}}

\newcommand\FF{{\mathcal F}}
\newcommand\GG{{\mathcal G}}
\newcommand\HH{{\mathcal H}}

\newcommand\LL{{\mathcal L}}
\newcommand\MM{{\mathcal M}}

\newcommand\PP{{\mathcal P}}

\newcommand\PMF{{\PP\kern-2pt\MM\FF}}

\newcommand\PML{{\PP\kern-2pt\MM\LL}}

\newcommand\hhat{\widehat}

\newcommand{\fsubd}{\mathrel{{\scriptstyle\searrow}\kern-1ex^d\kern0.5ex}}
\newcommand{\bsubd}{\mathrel{{\scriptstyle\swarrow}\kern-1.6ex^d\kern0.8ex}}
\newcommand{\fsubeq}{\mathrel{\raise-.7ex\hbox{$\overset{\searrow}{=}$}}}
\newcommand{\bsubeq}{\mathrel{\raise-.7ex\hbox{$\overset{\swarrow}{=}$}}}

\newcommand{\tsh}[1]{\left\{\kern-.9ex\left\{#1\right\}\kern-.9ex\right\}}

\begin{document}
\title{Relative Hyperbolic Extensions of Groups and Cannon-Thurston Maps}
\author{Abhijit Pal}
\maketitle
\begin{center}
{Stat-Math Unit, Indian Statistical Institute, 203 B.T.Road,\\ Kolkata 700108, Email: abhijit\_r@isical.ac.in}
\end{center}
\date{}
\begin{abstract}
Let $1\rightarrow (K,K_1)\rightarrow (G,N_G(K_1))\rightarrow(Q,Q_1)\rightarrow 1$ be a short exact sequence of pairs of finitely generated groups with $K$ strongly hyperbolic relative to proper subgroup $K_1$. Assuming  that for all 
$g\in G$ there exists $k\in K$ such that $gK_1g^{-1}=kK_1k^{-1}$, we prove that there exists a quasi-isometric section $s\colon Q \to G$. Further we prove that if $G$ is strongly hyperbolic relative to the normalizer subgroup $N_G(K_1)$ and weakly hyperbolic relative to $K_1$, then there exists a Cannon-Thurston map for the inclusion $i\colon\Gamma_K\to \Gamma_G$.
\end{abstract}

\section{Introduction}\label{sec:1}
Let us consider the short exact sequence of finitely generated groups 
\[
1\rightarrow K\rightarrow G\rightarrow Q\rightarrow 1.
\]
such that $K$ is non-elementary word hyperbolic. In \cite{mosher}, Mosher proved that if $G$ is hyperbolic, then $Q$ is hyperbolic.
To prove $Q$ hyperbolic, Mosher (in \cite{mosher}) constructed a quasi-isometric section  from $Q$ to $G$, that is, a map $s\colon Q\to G$ satisfying 
\[
\frac{1}{k}d_Q(q,q')-\epsilon \leq d_G(s(q),s(q'))\leq kd_Q(q,q')+\epsilon,
\]
for all $q,q'\in Q$, where $d_G$ and $d_Q$ are word metrics and $k\geq1$,$\epsilon\geq 0$ are constants. In \cite{mit0}, existence of a Cannon-Thurston map for the embedding $i\colon\Gamma_K\to\Gamma_G$ was proved, where $\Gamma_K$, $\Gamma_G$ are respectively the Cayley graphs of $K$, $G$.
Here in this paper, we will generalize these results to the case where the kernel is strongly hyperbolic relative to a cusp subgroup.

One of our main theorems states:\\ 
{\bf{Theorem} \ref{sec thm}}
Suppose we have a short exact sequence of finitely generated groups 
\[
1\rightarrow K\rightarrow G\stackrel{p}{\rightarrow} Q\rightarrow 1,
\]
with $K$ strongly  hyperbolic relative to a subgroup $K_1$ and for all $g\in G$ there exists $k\in K$ such that $gK_1g^{-1}=kK_1k^{-1}$ then there exists a $(k,\epsilon)$ quasi-isometric section $s\colon Q \to G$ for some constants $k\geq1$,$\epsilon\geq 0$.

As a corollary of this theorem, under the same hypothesis, we can take the image of quasi-isometric section to be in $N_G(K_1)$.

Let $S$ be a once-punctured torus then its fundamental group $\Pi_1(S)={\mathbb F}(a,b)$ is strongly hyperbolic relative to the peripheral subgroup $H=<aba^{-1}b^{-1}>$.
Let $M$ be a $3$-manifold fibering over circle with fiber  $S$ such that  the fundamental group $\Pi_1(M)$ is strongly hyperbolic relative to the subgroup $H\bigoplus\mathbb{Z}$. Then we have a short exact sequence of pairs of finitely generated groups:

\[
1\rightarrow(\Pi_1(S),H)\rightarrow(\Pi_1(M),H\bigoplus{\integers})\rightarrow(\integers,\integers)\rightarrow {1}.
\]
Let $K=\Pi_1(S)$, $G=\Pi_1(M)$ and let $\Gamma_K$, $\Gamma_G$ be Cayley graphs of $K$, $G$ respectively. 
Bowditch (in \cite{bow ct1}) and Mahan Mj. (in \cite{Mj bdd}), proved the existence of a Cannon-Thurston map for the embedding $i\colon\Gamma_K\to\Gamma_G$.

Motivated by this example, we will prove the following theorem :\\
\textbf{Theorem \ref{end thm}}
Suppose we have a short exact sequence of pairs of finitely generated groups 
\[
1\rightarrow (K,K_1)\stackrel{i}\rightarrow (G,N_G(K_1))\stackrel{p}{\rightarrow}(Q,Q_1)\rightarrow 1
\]
with $K$ strongly hyperbolic relative to the cusp subgroup $K_1$ and for all $g\in G$ there exists $k\in K$ such that $gK_1g^{-1}=kK_1k^{-1}$. If $G$ is  strongly  hyperbolic relative to $N_G(K_1)$ and weakly hyperbolic relative to the subgroup $K_1$, then there exists a Cannon-Thurston map for the embedding $i\colon \Gamma_K\to \Gamma_G$, where 
$\Gamma_K$ and $\Gamma_G$ are  Cayley graphs of $K$ and $G$ respectively.

\textbf{Acknowledgments:} I would like to thank my advisor, Mahan Mj., for suggesting me this problem and for his useful comments.

\section{ Relative Hyperbolicity and Quasi-isometric Section}\label{rel hyp}

For generalities of hyperbolic groups and hyperbolic metric spaces refer to \cite{hyp}.
Gromov defined  relatively hyperbolic group as follows:
\begin{defn}(Gromov \cite{gro})
Let $G$ be a finitely generated group acting freely and properly discontinously by isometries on a proper and $\delta$-hyperbolic metric space $X$,
such that the quotient space $X/G$ is quasi-isometric to $[0,\infty)$. Let $H$ denote the stabilizer subgroup of the endpoint on $\partial X$ of a lift of this ray to $X$. Then $G$ is said to be strongly hyperbolic relative to $H$. The subgroup $H$ is said to be a \textit{Parabolic or Cusp} subgroup and the end point on $\partial X$ as \textit{parabolic} end point.
\end{defn}

Thus for a group $G$ strongly hyperbolic relative to the subgroup $H$  there is  a natural bijective correspondence between parabolic end points and  parabolic subgroups of $G$. Infact, a parabolic end point corresponds to a subgroup of the form $aHa^{-1}$ for some $a\in G$. 

\begin{defn}(Farb \cite{farb}) Let $G$ be a finitely generated group, and let $H$ be a finitely generated subgroup of $G$. Let $\Gamma _G$ be the Cayley graph of $G$.
Let $\hhat \Gamma_G$ be a new graph obtained from $\Gamma_G$ as follows:

For each left coset $gH$ of $H$ in $G$, we add a new vertex $v(gH)$ to $\Gamma_G$, and add an edge $e(gh)$ of length $1/2$ from each element $gh$ of $gH$ to the vertex $v(gH)$. We call this new graph the \textit{Coned-off} Cayley graph of $G$ with respect to $H$, and 
denote it by $\hhat \Gamma_G=\hhat \Gamma_G(H)$.   

We say that $G$ is \textbf{weakly hyperbolic} relative to the subgroup $H$ if the Coned-off Cayley Graph $\hhat \Gamma_G$ is hyperbolic.
\end{defn}
Geodesics in the coned-off space $\hhat\Gamma_G$ will be called as electric geodesics.
For a path $\gamma \subset \Gamma_G$, there is an induced path $\hhat{\gamma}$ in $\hhat{\Gamma_G}$ obtained by identifying $\gamma$ with $\hhat{\gamma}$ as sets. 
If $\hhat{\gamma}$ is an electric geodesic (resp. $P$-quasigeodesic), $\gamma$ is called a {\em relative geodesic} (resp.{\em relative $P$-quasigeodesic}).
If $\hhat\gamma$ passes through some cone point $v(gH)$, we say that $\hhat\gamma$ {\em penetrates} the coset $gH$.
\begin{defn}(Farb \cite{farb})
  $\hhat\gamma$ is said to be an {\bf electric
$(K, \epsilon)$-quasigeodesic in (the electric space) $\hhat{\Gamma_G}$
 without backtracking } if
 $\hhat\gamma$ is an electric $K$-quasigeodesic in $\hhat{\Gamma_G}$ and
 $\hhat\gamma$ does not return to any left coset after leaving it.
\end{defn}

\begin{defn}(Farb \cite{farb}){\bf Bounded Coset Penetration Properties:}\label{bcp}
The pair $(G,H)$ is said to satisfy \textit{bounded coset penetration property} if, for every $P\geq 1$, there is a constant $D=D(P)$ such that if $\alpha$ and $\beta$ are two electric $P$-quasigeodesics without backtracking starting and ending at same points, then the following conditions hold :
\begin{enumerate}
\item If $\alpha$ penetrates a coset $gH$ but $\beta$ does not penetrate $gH$, then $\alpha$ travels a $\Gamma_G$-distance of at most $D$ in $gH$.
\item If both $\alpha$ and $\beta$ penetrate a coset $gH$, then  vertices in $\Gamma_G$ at which $\alpha$ and $\beta$ first enter $gH$ lie a $\Gamma_G$-distance of at most $D$ from each other; similarly for the last exit vertices.
\end{enumerate}
\end{defn}

Next we recall Farb's definition of relatively hyperbolic group (in the strong sense) from \cite{farb}:

\begin{defn}(Farb \cite{farb})
$G$ is said to be strongly  hyperbolic relative to $H$ if $G$ is weakly hyperbolic relative to $H$ and the pair $(G,H)$ satisfies bounded coset penetration property.
\end{defn}

%The phrase 'relatively hyperbolic' will mean 'strongly relatively hyperbolic' throughout this paper. 
\begin{thm}(Bowditch \cite{bow}, Szczepanski \cite{ski})
$G$ is strongly hyperbolic relative to the cusp subgroup $H$ in the sense of Farb if and only if $G$ is strongly hyperbolic relative to $H$ in the sense of Gromov.
\end{thm}

\begin{defn}(Bowditch {\cite{bow}})
{\bf Relative Hyperbolic Boundary:}
Let $G$ be a strongly  hyperbolic group relative to $H$, then by Gromov's definition $G$ acts properly discontinuously  on a proper hyperbolic space $X$. The relative hyperbolic boundary of $G$ is the Gromov boundary, $\partial X$, of $X$. We denote the relative hyperbolic boundary of the pair $(G,H)$ by $\partial \Gamma(G,H)$.
\end{defn}

Bowditch in \cite{bow} showed that the relative hyperbolic boundary $\partial \Gamma(G,H)$ is well defined.

\begin{defn}
Let $1\rightarrow K\rightarrow G\rightarrow Q\rightarrow 1$ be a short exact sequence of finitely generated groups with $K$ strongly  hyperbolic relative to $K_1$. We say that \textbf{$G$ preserves cusps} if $gK_1g^{-1}=a_gK_1a^{-1}_g$ for some $a_g\in K$. 
\end{defn}

\begin{defn}(Mosher \cite{mosher}) {\bf {Quasi-isometric section :}}
Let $1\rightarrow K\rightarrow G\rightarrow Q\rightarrow 1$ be a short exact sequence of finitely generated groups. A map $s\colon Q\to G$ is said to be a $(R,\epsilon)$ quasi-isometric section if 
\[
\frac{1}{R}d_Q(q,q')-\epsilon \leq d_G(s(q),s(q'))\leq Rd_Q(q,q')+\epsilon,
\]
for all $q,q'\in Q$, where $d_G$ and $d_Q$ are word metrics and $R\geq1$,$\epsilon\geq 0$ are constants.
\end{defn}

Let $K$ be a strongly hyperbolic group relative to a cusp subgroup $K_1$. For each parabolic point $\alpha\in \partial \Gamma(K,K_1)$, there is a unique subgroup of the form $aK_1a^{-1}$. Now Hausdorff distance between two sets $aK_1$ and $aK_1a^{-1}$ is uniformly bounded by the length of the word $a$. Hence $\alpha$ corresponds to a left coset $aK_1$ of $K_1$ in $K$.

Let $1\rightarrow K\stackrel{i}\rightarrow G\stackrel{p}{\rightarrow} Q\rightarrow 1$
be a short exact sequence of finitely generated groups with $K$  non-elementary and strongly hyperbolic relative to a subgroup $K_1$.

We use the following notations for our further purpose:
\begin{itemize}
\item Let $\Pi$ be the set of all parabolic end points for the  relatively hyperbolic group $K$ with cusp subgroup as $K_1$.

\item Let $\Pi^2= \{(\alpha_1,\alpha_2): \alpha_1~\mbox{and}~\alpha_2~\mbox{are distinct elements in}~\Pi \}$.

\item For $a\in K$, let $i_a\colon K\to K$ denote the inner automorphism and $L_a\colon K\to K$  the left translation.

\item For $g\in G$, let $I_g\colon K\to K$ be the outer automorphism, that is, $I_g(k)=gkg^{-1}$ and $L_g\colon G\to G$ be the left translation.
\end{itemize}

$G$ preserves cusps, so for each $g\in G$ there exists $a_g\in K$ such that $a_g^{-1}g\in N_G(K_1)$. If $b\in K$, then it can be easily proved that $d_K(a_gK_1,gbg^{-1}a_gK_1)\leq d_K(K_1,bK_1) + 2l_K(a_g^{-1}g)$. Since $I_g(bK_1)=g(bK_1)g^{-1}=gbg^{-1}a_gK_1a_g^{-1}$ and Hausdorff distance between $gbg^{-1}a_gK_1$ and $gbg^{-1}a_gK_1a_g^{-1}$ is bounded, $I_g$ will induce a map $\tilde I_g\colon \Pi \to \Pi$ and  $\tilde I_g$ is a bijection. Therefore, $\tilde I_g$ will induce a bijective map $\tilde I^2_g\colon \Pi^2\to \Pi^2$. For the sake of convenience of notation we will use $I_g$ for $\tilde I_g$ and $\tilde I^2_g$. Similarly, for $a\in K$, $i_a$ and $L_a$ will induce bijective maps (with same notation) from $\Pi$ to $\Pi$ and $\Pi^2$ to $\Pi^2$.

The following theorem is generalization of Mosher's technical result Quasi-isometric section lemma to the relatively hyperbolic case.

\begin{thm}\label{sec thm}
Suppose we have a short exact sequence of finitely generated groups 
\[
1\rightarrow K\stackrel{i}\rightarrow G\stackrel{p}{\rightarrow} Q\rightarrow 1,
\]
such that $K$ is  strongly hyperbolic relative to a subgroup $K_1$ and  $G$ preserves cusps, then there exists a  $(R,\epsilon)$ quasi-isometric section $s\colon Q \to G$ for some $R\geq1$,$\epsilon\geq 0$.
\end{thm}
\noindent{\bf Proof.} Let $\alpha=(\alpha_1,\alpha_2)\in \Pi^2$, then stabilizer subgroups of $\alpha_i$'s are $a_iK_1a_i^{-1}$ for some $a_i\in K$, where $i=1,2$. Let $\lambda$ be a relative geodesic in $\Gamma_K$ with starting at some point of $a_1K_1$ and ending at some point of $a_2K_1$ and $\hhat\lambda$ passing through cone points $v(a_iK_1)$, where $i=1,2$. Let $x$ be the exit point of $\lambda$ from the left coset $a_1K_1$. If $\mu$ is another such relative geodesic with end points same as $\lambda$ and $y$ as its exit point from $a_1K_1$, due to bounded coset penetration property 2, $d_{K}(x,y)\leq D$, where $D$ is the constant as in definition \ref{bcp}. Let $B_{\alpha}$ be the set of all exit points of relative geodesics $\lambda$ with starting at some point of $a_1K_1$ and ending at some point of $a_2K_1$ and $\hhat\lambda$ passing through $v(a_iK_1)$'s, $i=1,2$. Then $B_{\alpha}$ is a bounded set with diameter less than or equal to $D$.

Let $C=\{\alpha\in \Pi^2\colon e_K\in B_\alpha\}$, where $e_K$ is the identity element in $K$. We fix an element $\eta = (\eta_1,\eta_2)\in \Pi ^2$. Let $\Sigma=\{g\in G : \eta \in I_g(C)\}$. $\Sigma$ will be proved to be image of a quasi-isometric section.\\ We first prove that, for any $g\in G$, ${\cup}_{a\in K}I_{ga}(C)=\Pi^2$:

Let $\alpha = (\alpha_1,\alpha_2)\in \Pi^2$.  Now $\alpha_i$ corresponds to left coset $a_iK_1$, where $i=1,2$. Let $\lambda$ be a relative geodesic in $\Gamma_K$ with starting at some point of $a_1K_1$ and ending at some point of $a_2K_1$ and $\hhat\lambda$ passing through cone points $v(a_iK_1)$, $i=1,2$, and $x_\alpha$ its exit point from $a_1K_1$. Then $x_\alpha\in B_\alpha$. Now there exists $k\in K$ such that $L_k(x_\alpha)=e_K$. Since $L_k$ is an isometry $L_k(\lambda)$ will be a relative geodesic joining points from $ka_iK_1$'s, $i=1,2$, with $\hhat{L_k(\lambda)}$ containing cone points $v(ka_iK_1)$, $i=1,2$, and $e_K$ being the exit point of $L_k(\lambda)$ from $ka_1K_1$. There exists $\beta_i\in \Pi$ such that $\beta_i$ corresponds to left coset $ka_iK_1$, $i=1,2$. Then $\beta=(\beta_1,\beta_2)\in \Pi^2$ and $e_K\in B_\beta$. Therefore $L_k(\alpha)=\beta\in C$.
Since $L_k$ and $i_k$ are same on relative hyperbolic boundary, we have $i_k(\alpha)\in C$. Thus $\cup_{a\in K}(i_a(C))=\Pi^2$. Consequently, for any $g\in G$, $\cup_{a\in K}I_{ga}(C)=\cup_{a\in K}I_gi_a(C)=I_g(\cup_{a\in K}(i_a(C)))=I_g(\Pi^2)=\Pi^2$.

Now we prove that $p(\Sigma)=Q$:

Let $q\in Q$, then there exists $g\in G$ such that $p(g)=q$. Now  ${\cup}_{a\in K}I_{ga}(C)=\Pi^2$ for any $g\in G$. Therefore for $\eta\in \Pi^2$ there exists $a\in K$ such that $\eta\in I_{ga}(C)$. Hence $ga\in \Sigma$ and $p(ga)=p(g)=q$.

Now we prove that there exists constant $R\geq 1,\epsilon \geq 0$ such that for all $g,g'\in \Sigma$ 
\[
\frac{1}{R}d_Q(p(g),p(g'))-\epsilon \leq d_G(g,g') \leq R d_Q(p(g),p(g'))+\epsilon. 
\]

We can choose a finite symmetric generating set $S$ of $G$ such that $p(S)$ is also a generating set for $Q$. Obviously, $d_Q(p(g),p(g'))\leq d_G(g,g')$ for all $g,g'\in G$. To prove $d_G(g,g')\leq R d_Q(p(g),p(g'))+ \epsilon$ for all $g,g'\in \Sigma$, it suffices to prove that there exists $R\geq 1$ such that if $d_Q(p(g),p(g'))\leq 1$ for some $g,g'\in \Sigma$, then $d_G(g,g')\leq R$.

Let $d_Q(p(g),p(g'))\leq 1$ for some $g,g'\in \Sigma$. Then $g^{-1}g'=ka$ for some $k\in K$ and $a$ is either identity in $G$ or a generator of $G$. Since $g,g'\in \Sigma$, $I_g(C)\cap I_{g'}(C)\neq \Phi$. Hence $I_{ka}(C)\cap C = I_{g^{-1}g'}(C)\cap C \neq \Phi$. Now $I_{ka}=i_k(I_a)$, therefore $i_k(I_a(C))\cap C \neq \Phi$. 

For each $\alpha\in\Pi^2$, we choose an element $a_\alpha\in B_\alpha$. Define a map $F\colon\Pi^2\to\Gamma_K$ by $F(\alpha)=a_\alpha$.

Since $L_k$ is an isometry, for $k\in K$, $ka_{\alpha}\in B_{k\alpha}$ and hence
 \begin{equation}\label{hd}
d_K(a_{k\alpha},ka_\alpha)=d_K(F(k\alpha),kF(\alpha))\leq D,
\end{equation}
where $k\alpha$ denotes the image of $\alpha$ under the map $L_k\colon\Pi^2\to\Pi^2$.

Let $B_D(e_K)$ be the closed $D$-neighborhood of $e_K$. Now $F(C)$ is contained in union of $B_\alpha$'s containing identity $e_K$. Therefore $F(C)$ is contained in $B_D(e_K)$. Since $G$ preserves cusps, there exists $s\in K$ such that $F(I_a(C))$ is contained in union of $B_\alpha$'s containing $s$ and hence $F(I_a(C))\subset B_D(s)$, where $B_D(s)$ is a closed $D$-neighborhood of $s$.
From (\ref{hd}), Hausdorff distance between two sets $F(kI_a(C))$ and $kF(I_a(C))$ is bounded by $D$. For a set $A\subset\Gamma_K$, let $N_D(A)$ denotes the closed $D$-neighborhood of $A$. Thus \begin{eqnarray*}
F(kI_a(C))\subset N_D(kF(I_a(C)))=kN_D(F(I_a(C)))\subset kB_{2D}(s).
\end{eqnarray*}
Now $K$ acts properly discontinuously on $\Gamma_K$, therefore \begin{equation*}
B_D(e_K)\cap kB_{2D}(s)\neq \Phi                                                                          \end{equation*} for finitely many $k$'s in $K$. This implies $F(C)\cap F(kI_a(C))\neq\Phi$ for finitely many $k$'s in $K$. And hence $C\cap L_k(I_a(C))=C\cap kI_a(C)\neq\Phi$ for finitely many $k$'s in $K$. $L_k=i_k$ on relative hyperbolic boundary, so $C\cap (I_{ka}(C))\neq\Phi$ for finitely many $k$'s in $K$. Thus $g^{-1}g'= ka$ for finitely many $k$'s. Since number of generators of $G$ is finite, there exists a constant $R\geq 1$ such that $d_G(g,g')\leq R$. 

Now we define $s\colon Q\to G$ as follows:

Let $q\in Q$ and let there exists $g,g'\in \Sigma$ such that $p(g)=p(g')=q$. Then by above inequality $d_G(g,g')\leq R$. We choose one element $g\in p^{-1}(q)\cap\Sigma$ for each $q\in Q$ and define $s(q)=g$. 
Then  $s$ defines a single valued map satisfying :
\[
\frac{1}{R}d_Q(q,q')-\epsilon \leq d_G(s(q),s(q'))\leq Rd_Q(q,q')+\epsilon.
\]
for some constants $R\geq 1$, $\epsilon \geq 0$ and for all $q,q'\in Q$.\qed

\begin{cor}\label{pair thm}
Suppose we have a short exact sequence of pairs of finitely generated groups 

\[
1\rightarrow (K,K_1)\rightarrow (G,N_G(K_1))\stackrel{p}{\rightarrow}(Q,Q_1)\rightarrow 1
\]
with $K$  strong relative hyperbolic with respect to the cusp subgroup $K_1$. If $G$ preserves cusps, then $Q_1=Q$ and there is a quasi-isometric section $s\colon Q\to N_G(K_1)$ satisfying 

\[
\frac{1}{R}d_Q(q,q')-\epsilon \leq d_{N_G(K_1)}(s(q),s(q'))\leq Rd_Q(q,q')+\epsilon
\]
where $q,q'\in Q$ and $R\geq1$,$\epsilon \geq 0$ are constants.
Further, if $G$ is weakly hyperbolic relative to $K_1$, then $Q$ is hyperbolic.
\end{cor}
\noindent{\bf Proof.}Let $q\in Q$, then there exists $g\in G$ such that $p(g)=q$. Since $G$ preserves cusps, $gK_1g^{-1}=aK_1a^{-1}$ for some $a\in K$. Therefore $a^{-1}g\in N_G(K_1)$ and $q=p(a^{-1}g)\in Q_1$. Thus $Q_1=Q$.

Let $\Pi^2_{K_1}=\{(\alpha_1,\alpha_2)\in\Pi^2:\alpha_1~\mbox{corresponds to subgroup}~K_1\}$ and $C=\{\alpha\in\Pi^2_{K_1}:e_K\in B_\alpha\}$, where $B_\alpha$ is defined as in above theorem. We fix an element $\eta\in\Pi^2_{K_1}$ and set $\Sigma=\{g\in N_G(K_1): \eta\in I_g(C)$. We choose a finite generating set $S$ of $G$ such that it contains a finite generating set of $N_G(K_1)$ and $p(S)$ is also a generating set of $Q$. Using  argument same as above theorem, by replacing $G$ with $N_G(K_1)$, we get a quasi-isometric section 
 $s\colon Q\to N_G(K_1)$ satisfying :
\[
\frac{1}{R}d_Q(q,q')-\epsilon \leq d_{N_G(K_1)}(s(q),s(q'))\leq Rd_Q(q,q')+\epsilon.
\]
for some constants $R\geq 1$, $\epsilon \geq 0$ and for all $q,q'\in Q$.

Since $d_Q(q,q')\leq d_G(s(q),s(q'))\leq d_{N_G(K_1)}(s(q),s(q'))$, we can take the quasi-isometric section $s\colon Q \to N_G(K_1)$ such that 

$\frac{1}{R}d_Q(q,q')-\epsilon \leq d_G(s(q),s(q'))\leq Rd_Q(q,q')+\epsilon$.

Now, let $\Gamma_G^{pel}$ denotes the electrocuted space obtained from $\Gamma_G$ by coning left cosets of $K_1$ in $G$. Since $G$ is weakly hyperbolic with respect to $K_1$, $\Gamma_G^{pel}$ is hyperbolic. We will prove that $Q$ is hyperbolic.

The quasi-isometric section $s\colon Q\to N_G(K_1)(\subset G)$ will induce a map 
$\hhat s\colon Q\to \Gamma_G^{pel}$. Now for all $q,q'\in Q$, 
$d_{G^{pel}}(\hhat s(q),\hhat s(q'))\leq d_G(s(q),s(q'))\leq R~d_Q(q,q')+\epsilon$, where $d_{G^{pel}}$ is the metric on $\Gamma_G^{pel}$. Obviously, $d_Q(q,q')\leq d_{G^{pel}}(\hhat s(q),\hhat s(q'))$. Hence $\hhat s$ is a quasi-isometric section from $Q$ to $\Gamma_G^{pel}$. Therefore $s(Q)$ is quasiconvex in $\Gamma_G^{pel}$. Since $\Gamma_G^{pel}$ is hyperbolic, $Q$ is hyperbolic.\qed

\section{Existence of Cannon-Thurston Maps}\label{c-t map}

Let $G$ be a group strongly hyperbolic relative to a subgroup $G_1$. Then there exists a complete hyperbolic space $Z$ such that $G$ acts properly discontinuously on $Z$ by isometries and $Z/G$ is quasi-isometric to $[0,\infty)$. Let $p$ be the parabolic end point corresponding to a lift of this ray. Then the subgroup stabilising $p$ is equal to $gG_1g^{-1}$ for some $g\in G$. Let $\Pi$ be the set of all parabolic end points and for each $p\in \Pi$, let $B(p)$ be a closed $G_1$ invariant such that $B(p)\cap B(q)=\Phi$ for all distinct parabolic end points $p,q$. Let $X'=Z\setminus\bigcup_{p\in\Pi}B(p)$. Then $X'$ is quasi-isometric to $\Gamma_G$. 

Let $X=\Gamma_G$ and $H_{gG_1}$ be the closed set in $\Gamma_G$ corresponding to the left coset $gG_1$ of $G_1$ in $G$. Let $\mathcal{H}_G$ $=\{H_{gG_1}:g\in G\}$. Let $X_h$ be the space obtained from $X$ by gluing $H_{gG_1}\times [0,\infty)$ to $H_{gG_1}$ for all $H_{gG_1}\in \mathcal{H}_G$ where  $H_{gG_1}\times [0,\infty)$ is equipped with the path metric $d_h$ such that
\begin{enumerate}
\item $d_{h,t}((x,t),(y,t)) = 2^{-t}d_H(x,y)$, where $d_{h,t}$ is the induced path
metric on $H_t=H\times \{t\}$.  

\item $d_h((x,t),(x,s))=\vert t-s \vert$ for all $x\in H$ and for all $t,s\in [0,\infty)$.
\end{enumerate}
Then $X_h$ is quasi-isometric to $Z$.

Thus if $G$ is strongly hyperbolic relative to $G_1$,  there exists a complete hyperbolic metric space $X_h$ such that $G$ acts properly discontinuously on $X_h$ and $X_h$ is obtained from $X$ as above. $H_{gG_1}\in \mathcal{H_G}$ was referred to as \textbf{horosphere-like sets} by Mahan Mj. and Lawerence Reeves in \cite{Mj com} and $H_{gG_1}\times [0,\infty)$ was referred to as \textbf{hyperbolic cones or horoball-like sets} in \cite{pal}.

\begin{defn}\label{Mj ct}
Let $\hat \lambda$ be an electric quasi-geodesic in the electric space $\hat X$ without backtracking.   For any  $H$ penetrated by $\hat \lambda$,  let $x_H$ and $y_H$ be the first entry
point  and the last exit point of $\hat \lambda$.   We join $x_H$ and $y_H$ by a hyperbolic geodesic segment in $H_h$ (identifying $\hhat{X}$ with the space obtained from $X_h$ by coning off the $H_h$'s.   This results in a path $\lambda_h$ in $X_h$. The path $\lambda$ will be called an {\bf electro-ambient quasigeodesic}.  
\end{defn}

\begin{lem}
An electro-ambient quasigeodesic is a quasigeodesic in the hyperbolic space $X_h$.  
\end{lem}

Consider the inclusion between pairs of relatively hyperbolic group $(H,H_1)\stackrel{i}\hookrightarrow (G,G_1)$. $i$ will induce a proper embedding $i\colon \Gamma_H \to \Gamma_G$. Let $X=\Gamma_G$ and $Y=\Gamma_H$. Recall that $X_h$ is the space obtained from $X$ by gluing the hyperbolic cones. 
Inclusion of a horosphere-like set in its hyperbolic cone is  uniformly proper, therefore inclusion of $X$ in $X_h$ is uniformly proper i.e. for all $M>0$ and $x,y\in X$ there exists $N>0$ such that if $d_{X_h}(x,y)\leq M$ then $d_G(x,y)=d_X(x,y)\leq N$.

Since $G$ preserves cusps, $i$ will induce a proper embedding $i_h\colon Y_h\to X_h$.  

\begin{defn}\label{cann-thu map}
A Cannon-Thurston map for $i\colon (\Gamma_H,{\mathcal{H_H}})\to (\Gamma_G,\mathcal{H_G})$ is said to exist if there exists a continuous extension $\tilde{i_h}\colon Y_h\cup \partial Y_h\to X_h\cup \partial X_h$ of $i_h\colon Y_h\to X_h$.
\end{defn}

To prove the existence of Cannon-Thurston map for the inclusion $i\colon (K,K_1)\to (G,N_G(K_1))$, we need the notion of Partial Electrocution.
\begin{defn}[Mahan Mj. and Lawrence Reeves \cite{Mj com}](Partial Electrocution)
Let $(X, \HH , \GG , \LL )$ be an ordered quadruple such that the
following holds:
\begin{enumerate}
\item $X$ is  hyperbolic relative to a collection of subsets
$H_\alpha$. 
\item For each $H_\alpha$ there is a uniform large-scale
retraction $g_\alpha : H_\alpha \rightarrow L_\alpha$ to some
(uniformly) $\delta$-hyperbolic metric space $L_\alpha$, i.e. there
exist $\delta , K, \epsilon > 0$ such that for all $H_\alpha$ there exists
a $\delta$-hyperbolic $L_\alpha$ and a map 
$g_\alpha : H_\alpha \rightarrow L_\alpha$ with
$d_{L_\alpha} (g_\alpha (x), g_\alpha (y)) \leq Kd_{H_\alpha}(x,y)
+ \epsilon $ for all $x, y \in H_\alpha$. Further, we denote the
collection of such $g_\alpha$'s as $\GG$.  
\end{enumerate}
 The {\bf partially electrocuted space} or
{\em partially coned off space} corresponding to $(X, \HH , \GG , \LL)$ 
is 
obtained from $X$ by gluing in the (metric)
mapping cylinders for the maps
 $g_\alpha : H_\alpha \rightarrow L_\alpha$.
\end{defn}
\begin{lem}(Mahan Mj. and Lawrence Reeves \cite{Mj com})\label{pel-track}
Given $K, \epsilon \geq 0$, there exists $C > 0$ such that the following
holds: \\
Let $\gamma_{pel}$ and $\gamma$ denote respectively a $(K, \epsilon )$
partially electrocuted quasigeodesic in $(X,d_{pel})$ and a hyperbolic
$(K, \epsilon )$-quasigeodesic in $(X_h,d_h)$ joining $a, b$. Then $\gamma \cap X$
lies in a (hyperbolic) $C$-neighborhood of (any representative of) 
$\gamma_{pel}$. Further, outside of  a $C$-neighborhood of the horoballs
that $\gamma$ meets, $\gamma$ and $\gamma_{pel}$ track each other.
\end{lem}

Let $G$ be a group strongly hyperbolic relative to the subgroup $G_1$.
Let $X_h$ be the complete hyperbolic space obtained from $X$.  We describe a special type of quasigeodesic in $X_h$ which will be essential for our purpose:
\begin{defn}(Mahan Mj. \cite{Mj ct})
We start with an electric quasi-geodesic $\hat \lambda$ in the electric space $\hat X$ without backtracking.   For any horosphere-like set $H$ penetrated by $\hat \lambda$,  let $x_H$ and $y_H$ be the respective entry and exit points to $H$.   We join $x_H$ and $y_H$ by a hyperbolic geodesic segment in $H\times [0, \infty)$.   This results a path in,  say $\lambda$,  in $X_h$.   The path $\lambda$ will be called an {\bf electro-ambient quasigeodesic}.  
\end{defn}

\begin{thm}(Mahan Mj. \cite{Mj ct})\label{e-a thm}:
An electro-ambient quasigeodesic is a quasigeodesic in the hyperbolic space $X_h$.  

\end {thm}
For the rest of paper, we will work with the following pair of short exact sequence of finitely generated groups :

\[
1\rightarrow (K,K_1)\rightarrow (G,N_G(K_1))\stackrel{p}{\rightarrow}(Q,Q_1)\rightarrow 1
\]

Since all groups are finitely generated, we can choose a finite generating set $S$ of $G$ such that $S$ contains finite generating set of $K,K_1,N_G(K_1)$ and $p(S)$ is also a finite generating set of $Q$.
We will assume the hypothesis of Theorem \ref{pair thm}. As a consequence, $Q_1=Q$ and there exists a $(R,\epsilon)$ quasi-isometric section $s\colon Q \to N_G(K_1)$ such that 

$\frac{1}{R}d_Q(q,q')-\epsilon \leq d_G(s(q),s(q'))\leq Rd_Q(q,q')+\epsilon$.

Further, we assume that $G$ is strongly hyperbolic relative to the subgroup $N_G(K_1)$.

Also, using a left translation $L_k$ by an element $k\in K_1$, we can assume that $s(Q)$ contains identity element $e_K$ of $K$ and $s(Q)\subset N_G(K_1)$.

We have assumed that $G$ is weakly hyperbolic relative to the subgroup $K_1$ and hence  the coned-off space $\Gamma_G^{pel}$ obtained by coning left cosets $gK_1$ of $K_1$  to a point $v(gK_1)$ is hyperbolic. As $\Gamma_Q$ is quasi-isometrically embedded in $\Gamma_G^{pel}$, $Q$ is hyperbolic. We have also assumed that $G$ is strongly hyperbolic relative to $N_G(K_1)$. Thus $\Gamma_G^{pel}$ becomes a partially electrocuted space obtained from $\Gamma_G$ by partial electrocuting the closed sets (horosphere-like sets) $H_{gN_G(K_1)}$ in $\Gamma_G$ corresponding to the left cosets $gN_G(K_1)$ to the hyperbolic space $g(s(Q))$, where $g(s(Q))$ denotes the image of $s(Q)$ under the left translation $L_g$ for $g\in G$.

Let $\lambda^b=\hhat\lambda\setminus\mathcal{H_K}$ denotes the portions of $\hhat\lambda$ that does not penetrate  horosphere-like sets in $\mathcal{H_K}$. The following Lemma gives a sufficient condition for the existence of Cannon-Thurston map for the inclusion $i\colon (\Gamma_K,{\mathcal{H_K}})\to (\Gamma_G,\mathcal{H_G})$. For instance see \cite{pal}.

\begin{lem}\label{suff cond}\cite{pal}
A Cannon-Thurston map for $i\colon (\Gamma_K,{\mathcal{H_K}})\to (\Gamma_G,\mathcal{H_G})$ exists if there exists a non-negative function $M(N)$ with $M(N)\rightarrow\infty$ as $N\rightarrow\infty$ such that the following holds:

Given $y_0\in \Gamma_K$ and an electric quasigeodesic segment $\hhat \lambda$ in $\hhat \Gamma_K$ if $\lambda^b=\hhat\lambda\setminus\mathcal{H_K}$  lies outside an  $N$-ball around $y_0\in \Gamma_K$, then for any partially electrocuted quasigeodesic $\beta_{pel}$ in $\Gamma_G^{pel}$ joining end points of $\hhat \lambda$, $\beta^b=\beta_{pel}\setminus\mathcal{H_G}$ lies outside an $M(N)$-ball around $i(y_0)$ in $\Gamma_G$.
\end{lem}

\subsection{Construction of Quasiconvex Sets and Retraction Map}
\label{Qc set}
Recall that for $g\in G$, $L_g\colon G\to G$ denotes the left translation by $g$ and $I_g\colon K\to K$ denotes the outer automorphism. Let $\phi_g= I_{g^{-1}}$ then $\phi_g(a)=g^{-1}ag$. Since $L_g$ is an isometry, $L_g$ preserves distance between left cosets of $G_1$ in $G$. Hence $L_g$ induces a isometry 
$\hhat L_g\colon \Gamma_G^{pel}\to\Gamma_G^{pel}$. The embedding $i\colon\Gamma_K\to \Gamma_G$ will induce an embedding $\hat i\colon\hhat\Gamma_K\to\Gamma_G^{pel}$.

Let $\hhat\lambda$ be an electric geodesic segment in $\hhat \Gamma_K$ with end points $a$ and $b$ in $\Gamma_K$. Let $\hhat \lambda_g$ be an electric geodesic in $\hhat \Gamma_K$ joining  $\phi_g(a)$ and $\phi_g(b)$.

Define
\[
 B_{\hhat\lambda}=\bigcup_{g\in s(Q)}{\hhat L_g.\hat{i}(\hhat \lambda_g)}.
\]

On $\hhat\Gamma_K$, define a map $\pi_{\hhat\lambda_g}\colon\hhat\Gamma_K\to \hhat\lambda_g$ taking $k\in \hhat\Gamma_K$ to one of the points on $\hhat\lambda_g$ closest to $k$ in the metric $d_{\hhat K}$.

\begin{lem}(Mitra \cite{mit1})\label{npp}
For $\pi_{\hhat\lambda_g}$ defined above, 
\[
d_{\hhat K}(\pi_{\hhat\lambda_g}(k),\pi_{\hhat\lambda_g}(k'))\leq Cd_{\hhat K}(k,k')+C
\]
for all $k,k'\in \hhat\Gamma_K$, where $C$ depends only on the hyperbolic constant of $\hhat\Gamma_K$.
\end{lem}

Now define $\Pi_{\hhat\lambda}\colon \Gamma_G^{pel}\to B_{\hhat\lambda}$ as follows:

For $g\in s(Q)$, $\Pi_{\hhat\lambda}.{\hhat L_g}.{\hhat i}(k)= \hhat {L_g}.\hhat {i}.\pi_{\hhat\lambda_g}(k)$ for all $k\in \Gamma_{\hhat K}$.

For every $g'\in \Gamma_G$, there exists a unique $g\in s(Q)$ such that $g'=L_g( i(k))$ for some unique $k\in K$. Hence, $\Pi_{\hhat\lambda}$ is well defined on the entire  space $\Gamma_G^{pel}$.

\begin{thm}(Mitra \cite{mit0,mit1})\label{ret}
There exists $C_0>0$ such that
\[
d_{\hhat G}(\Pi_{\hhat\lambda}(g),\Pi_{\hhat\lambda}(g'))\leq C_0d_{\hhat G}(g,g')+C_0
\]
for all $g,g'\in \Gamma_G^{pel}$. In particular, if $\Gamma_G^{pel}$ is hyperbolic then $B_{\hhat\lambda}$ is uniformly (independent of $\hhat\lambda$) quasiconvex.
\end{thm}

\subsection{Proof of Theorem}\label{main thm}
Since $i\colon \Gamma_K\to\Gamma_G$ is an embedding we identify $k\in K$ with its image $i(k)$.

Let \begin{itemize}
\item $\hhat \mu_g=\hhat L_g(\hhat\lambda_g)$, where $g\in s(Q)$.
\item $\mu^b_g=\hhat\mu_g\setminus\mathcal{H_G}$.
\item $B_{\lambda^b}=\bigcup_{g\in s(Q)}\mu^b_g$.
\item $Y=\Gamma_K$ and $X=\Gamma_G$.
\end{itemize}
\begin{lem}\label{end lem}
There exists $A>0$ such that for all $x\in \mu^b_g\subset B_{\lambda^b}\subset B_{\hhat\lambda}$ if $\lambda^b$ lies outside $B_N(p)$ for a fixed reference point $p\in\Gamma_K$ then $x$ lies outside an $\frac{f(N)}{A+1}$ ball about $p$ in $\Gamma_G$, where $f(N)\rightarrow\infty$ as $N\rightarrow\infty$.
\end{lem}

\noindent {\bf Proof.} Let $x\in \mu^b_g$ for some $g\in s(Q)$. Let $\gamma$ be a geodesic path in $\Gamma_Q$ joining the identity element $e_Q$ of $\Gamma_Q$ and $p(x)\in\Gamma_Q$. Order the vertices on $\gamma$ so that we have a finite sequence $e_Q=q_0,q_1,...,q_n=p(x)=p(g)$ such that $d_Q(q_i,q_{i+1})=1$ and $d_Q(e_Q,p(x))=n$. Since $s$ is a quasi-isometric section, this gives a sequence $s(q_i)=g_i$ such that $d_G(g_i,g_{i+1})\leq R+\epsilon = R_1\mbox{ (say)}$. Observe that $g_n=g$ and $g_0=e_G$. Let $B_{R_1}(e_G)$ be a closed ball around $e_G$ of radius $R_1$, then  $B_{R_1}(e_G)$ is finite. Now for each $g\in G$, the outer automorphism $\phi_g$ is a quasi-isometry. Thus there exists $K\geq 1$ and $\epsilon\geq 0$ such that for all $g\in  B_{R_1}(e_G)$, $\phi_g$ is a $(K,\epsilon)$ quasi-isometry and $K,\epsilon$ are independent of elements of $G$. Let $s_i=g^{-1}_{i+1}g_i$, then $s_i\in B_{R_1}(e_G)$, where $i=0,...,n-1$. Hence $\phi_{s_i}$ is $(K,\epsilon)$ quasi-isometry.

Since $s(Q)\subset N_G(K_1)$, we have $s_i\in N_G(K_1)$ for all $i$. Therefore $\phi_{s_i}$ will induce a $(\hat K,\hat\epsilon)$ quasi-isometry $\hhat\phi_{s_i}$ from $\hhat \Gamma_K$ to $\hhat \Gamma_K$, where $\hat K,\hat \epsilon$ depends only $K$ and $\epsilon$.

Now $x\in  \mu^b_{g_n}$ and $L_g$ preserves distance between left cosets for all $g\in G$, hence there exists $x_1\in \lambda^b_{g_n}$ such that $x=L_{g_n}(x_1)$.

 %$L_{g_{n-1}}\phi_{s_{n-1}}(x_1)=g_{n-1}(g_n^{-1}g_{n-1})^{-1}x_1(g_n^{-1}g_{n-1})=xs_{n-1}$

Let $[p,q]_{g_n}\subset \lambda^b_{g_n}$ be the connected portion of $\lambda^b_{g_n}$ on which $x_1$ lies. Since $\hhat\phi_{s_{n-1}}$ is a quasi-isometry,  $\hhat \phi_{s_{n-1}}([p,q]_{g_n})$ will be an electric quasigeodesic lying outside horosphere-like sets and hence it is a quasigeodesic in $Y_h$ lying at a uniformly bounded distance $\leq C_1$  from $\lambda^h_{g_{n-1}}$ in $Y_h$ (and hence in $X_h$), where $\lambda^h_{g_{n-1}}$ is electroambient representative of $\hhat\lambda_{g_{n-1}}$ and $Y_h, X_h$ are respectively  the complete hyperbolic metric spaces corresponding to $Y,X$. Thus there exist $x_2\in \lambda^h_{g_{n-1}}$ such that $d_{X_h}(\phi_{s_{n-1}}(x_1),x_2)\leq C_1$. But $x_2$ may lie inside horoball-like set penetrated by $\hhat\lambda_{g_{n-1}}$. Due to bounded coset (horosphere) penetration properties there exists $y\in\lambda^b_{g_{n-1}}$ such that $d_{X_h}(x_2,y)\leq D$.

Thus $d_{X_h}(\phi_{s_{n-1}}(x_1),y)\leq C_1+D$. Since $X=\Gamma_G$ is properly embedded in $X_h$, there exists $M>0$ depending only upon $C_1,D$ such that $d_G(\phi_{s_{n-1}}(x_1),y)\leq M$.

Hence $d_G(L_{g_{n-1}}(\phi_{s_{n-1}}(x_1)),L_{g_{n-1}}(y))=d_G(\phi_{s_{n-1}}(x_1),y)\leq M$ and $L_{g_{n-1}}(y)\in \mu^b_{g_{n-1}}$.

Let $z=L_{g_{n-1}}(y)$, then
\begin{eqnarray*}
d_G(x,z)&\leq& d_G(x,L_{g_{n-1}}(\phi_{s_{n-1}}(x_1)))+d_G(L_{g_{n-1}}(\phi_{s_{n-1}}(x_1)),L_{g_{n-1}}(y))\\&\leq& d_G(x,xs_{n-1})+M\\&\leq& R_1+M=A (say). 
\end{eqnarray*}
Thus, we have shown that for $x\in\mu^b_{g_n}$ there exists $z\in\mu^b_{g_{n-1}}$ such that $d_G(x,z)\leq A$. Proceeding in this way, for each $y\in\mu^b_{g_i}$ there exists $y'\in\mu^b_{g_{i-1}}$ such that $d_G(y,y')\leq A$.

Hence there exists $x'\in\lambda^b$ such that $d_G(x,x')\leq An$.

Since $\Gamma_K$ is properly embedded in $\Gamma_G$ there exists $f(N)$ such that $\lambda^b$ lies outside $f(N)$-ball about $p$ in $\Gamma_G$ and $f(N)\rightarrow\infty$ as $N\rightarrow\infty$.

Therefore $d_G(x',p)\geq f(N)$.

Thus
 \begin{eqnarray*}
    d_G(x,p)\geq f(N)-d_G(x,x')\geq f(N)-An.
\end{eqnarray*}

Also $d_G(x,p)\geq n$.

Therefore, $d_G(x,p)\geq \frac{f(N)}{A+1}$, that is, $x$ lies outside $\frac{f(N)}{A+1}$-ball about $p$ in $\Gamma_G$.\qed

\begin{thm}\label{end thm}
Given a short exact sequence of pairs of finitely generated groups 

\[
1\rightarrow (K,K_1)\rightarrow (G,N_G(K_1))\stackrel{p}{\rightarrow}(Q,Q_1)\rightarrow 1
\]
with $K$ non-elementary and strongly hyperbolic relative to the cusp subgroup $K_1$. If $G$ preserves cusps,  strongly  hyperbolic relative to $N_G(K_1)$ and weakly hyperbolic relative to the subgroup $K_1$, then there exists a Cannon-Thurston map for the embedding $i\colon \Gamma_K\to \Gamma_G$, where 
$\Gamma_K$ and $\Gamma_G$ are  Cayley graphs of $K$ and $G$ respectively.
\end{thm}
\noindent {\bf Proof.} It suffices to prove the condition of Lemma \ref{suff cond}.

So for a fixed reference point $p\in \Gamma_K$, we assume that $\hhat\lambda$ is an electric geodesic segment in $\hhat\Gamma_K$ such that $\lambda^b(\subset \Gamma_K)$ lies outside an $N$-ball $B_N(p)$ around $p$. Let $\beta_{pel}$ be a quasigeodesic in the partially electrocuted space $\Gamma_G^{pel}$ joining the end points of $\hhat\lambda$. Recall from Theorem \ref{ret} that $\Pi_{\hhat\lambda}$ is a retraction map from $\Gamma_G^{pel}$ onto the quasiconvex set $B_{\hhat\lambda}$ which satisfies  the Lipschitz's condition. Let $\beta'_{pel}=\Pi_{\hhat\lambda}(\beta_{pel})$, then $\beta'_{pel}$ is a quasigeodesic in $\Gamma_G^{pel}$ lying on $B_{\hhat\lambda}$. So $\beta'_{pel}$ lies in a $P$-neighborhood of $\beta_{pel}$ in $\Gamma_G^{pel}$. But $\beta'_{pel}$ might backtrack. $\beta'_{pel}$ can be modified to form a quasigeodesic $\gamma_{pel}$ in $\Gamma_G^{pel}$ of the same type (i.e. lying in a $P$-neighborhood of $\beta_{pel}$) without backtracking with end points remaining the same.
By Lemma \ref{pel-track}, $\beta_{pel}$ and $\gamma_{pel}$ satisfy bounded coset (horosphere) penetration  properties with the closed sets (horosphere-like sets) in $\Gamma_G$  corresponding to the left cosets  of $N_G(K_1)$ in $G$. Thus if $\gamma_{pel}$ penetrates a horosphere-like set $C_{gN_G{(K_1)}}$ corresponding to the left coset $gN_G(K_1)$ of $N_G(K_1)$ in $G$ and $\beta _{pel}$ does not, then  length of geodesic traversed by $\gamma^h$, where $\gamma^h$ is the electroambient quasigeodesic representative of $\gamma_{pel}$, inside $C_{gN_G(K_1)}\times [0,\infty)$ is uniformly bounded. Let ${\mathcal{C}}=\{C_{gN_G{(K_1)}}:g\in G\}$.
Thus there exists $C_1\geq 0$ such that if $x\in \beta ^b_{pel} = \beta _{pel}\setminus \mathcal C$,
then there exists $y\in \gamma ^b_{pel} = \gamma _{pel}\setminus \mathcal C$ such that 
$d_G(x,y)\leq C_1$.

Since $y\in\gamma ^b_{pel}\subset B_{\lambda^b}$, by Lemma \ref{end lem}, $d_G(y,p)\geq \frac{f(N)}{A+1}$.

Therefore, $d_G(x,p)\geq  \frac{f(N)}{A+1}-C_1 (=M(N),\mbox{ say})$ and $M(N)\rightarrow\infty$ as $N\rightarrow\infty$. By Lemma \ref{suff cond}, a Cannon-Thurston map for $i\colon\Gamma_K\to\Gamma_G$ exists.\qed

\end{document}